\documentclass[reqno]{amsart}
\usepackage{amsfonts}
\usepackage{amsmath}
\usepackage{amssymb}
\usepackage{mathrsfs}
\usepackage[colorlinks]{hyperref}
\usepackage{tikz}
\usetikzlibrary{matrix,shadows,positioning}
\usepackage{algorithm}
\usepackage{algpseudocode}
\usepackage{algorithmicx}
\usepackage{graphicx,lscape}
\usepackage{subcaption}
\DeclareGraphicsExtensions{.pdf,.png,.jpg}

\newcommand{\Lcal}{\mathcal L}
\newcommand{\Ah}{A_h}
\newcommand{\ip}[2]{\left\langle #1,#2\right\rangle_h}
\newcommand{\normh}[1]{\left\|#1\right\|_h}
\newcommand{\normAh}[1]{\left\|#1\right\|_{A_h,h}}

\begin{document}
\title[Identification of a time-dependent coefficient]
{Numerical identification of a time-dependent coefficient in a heat equation with the spectral fractional Laplacian}

\author[Arshyn Altybay]{Arshyn Altybay}
\address{
Arshyn Altybay (corresponding author):
\endgraf
Institute of Mathematics and Mathematical Modeling,
\endgraf
125 Pushkin str., 050010 Almaty, Kazakhstan
\endgraf
and
\endgraf
International Engineering and Technological University, Almaty, Kazakhstan
\endgraf
{\it E-mail address} {\rm arshyn.altybay@gmail.com, arshyn.altybay@math.kz}
}

\author[G. Nalzhupbayeva]{Gulzat Nalzhupbayeva}
\address{
Gulzat Nalzhupbayeva:
\endgraf
al--Farabi Kazakh National University
\endgraf
71 al--Farabi ave., 050040 Almaty, Kazakhstan
\endgraf
and
\endgraf
Institute of Mathematics and Mathematical Modeling
\endgraf
125 Pushkin str., 050010 Almaty, Kazakhstan
\endgraf
{\it E-mail address:} {\rm gulzat.nalzhupbayeva@gmail.com}
}

\thanks{}
\subjclass[2020]{65M06, 65M12, 35R30, 35K20}
\keywords{inverse coefficient problem, spectral fractional Laplacian, fractional heat equation, integral overdetermination, Crank--Nicolson method, fractional matrix power}

\begin{abstract}
We consider the inverse problem of identifying a time-dependent source coefficient in a one-dimensional heat equation governed by the spectral Dirichlet fractional Laplacian from a weighted integral measurement. We present the continuous semigroup formulation and reduce the inverse
problem to a Volterra equation under a nondegeneracy condition. For the numerical approximation, the spectral fractional operator is discretised by the fractional matrix power $A_h=L_h^s$ of the standard Dirichlet difference Laplacian and combined with a Crank--Nicolson scheme. We prove unconditional stability and an $O(\tau^2+h^2)$ convergence estimate under a second-order spatial consistency assumption. A scalar reconstruction formula for the unknown coefficient is derived together with a discrete identifiability condition, and conditional convergence of the coupled state--coefficient reconstruction is established. For noisy integral measurements, regularised differentiation is used to approximate the required derivative data. Numerical experiments
illustrate stable reconstruction for \(1\%\)--\(5\%\) relative noise and are consistent with the derived conditional error estimate.
\end{abstract}

\maketitle
\numberwithin{equation}{section}
\newtheorem{theorem}{Theorem}[section]
\newtheorem{lemma}[theorem]{Lemma}
\newtheorem{corollary}[theorem]{Corollary}
\newtheorem{remark}[theorem]{Remark}
\newtheorem{definition}[theorem]{Definition}
\newtheorem{problem}[theorem]{Problem}
\newtheorem{proposition}[theorem]{Proposition}
\newtheorem{assumption}[theorem]{Assumption}
\allowdisplaybreaks

\section{Introduction}\label{1}

Fractional diffusion equations provide a useful mathematical framework for
describing spatially nonlocal transport and anomalous diffusion phenomena \cite{Kwasnicki17,Bouchaud90, Valdinoci09,Daoud22}.
On bounded domains, models involving fractional powers of elliptic operators
are particularly attractive because they retain the geometry and boundary
conditions of the underlying classical problem while incorporating nonlocal
spatial effects. In many applications, however, the governing equation
contains parameters or source intensities that cannot be measured directly
and must instead be identified from indirect observations. This leads
naturally to inverse coefficient and inverse source problems for fractional
diffusion equations.

A practically important situation occurs when the available observation is
not the complete state $u(t,x)$ at every spatial point, but only an aggregate
or weighted measurement of the form
\[
\int_{\Omega}u(t,x)\omega(x)\,dx.
\]
Such measurements arise naturally when pointwise observations are
unavailable, unreliable, or replaced by spatially averaged sensor data.
Recovering a time-dependent source intensity from this type of information
is therefore both mathematically challenging and relevant for parameter
identification in distributed systems.

Let $\Omega=(0,l)$ and $T>0$. In this work, we consider the simultaneous
determination of a time-dependent source coefficient $r(t)$ and the state
$u(t,x)$ satisfying
\begin{equation}\label{1.1}
 u_t(t,x)+(-\Delta_D)^s u(t,x)=r(t)f(t,x),
 \qquad (t,x)\in Q:=(0,T]\times\Omega,
\end{equation}
subject to the initial condition
\begin{equation}\label{1.2}
 u(0,x)=\varphi(x),\qquad x\in\Omega,
\end{equation}
and the homogeneous Dirichlet boundary condition
\begin{equation}\label{1.3}
 u(t,0)=u(t,l)=0,\qquad 0<t\le T.
\end{equation}
Here the homogeneous boundary condition is understood through the Dirichlet
realisation of $-\Delta$; whenever the corresponding traces are defined, it
coincides with the classical zero boundary condition. The additional
information is given by the weighted integral measurement
\begin{equation}\label{1.4}
 \int_0^l u(t,x)\omega(x)\,dx=w(t),
 \qquad 0\le t\le T.
\end{equation}
Here $0<s<1$, the function $f(t,x)$ is known, $\varphi$ is the prescribed
initial state, $\omega$ is a given observation weight, and $w(t)$ denotes the
measured response. The inverse problem is to recover the pair $(u,r)$ from
the data $\{f,\varphi,\omega,w\}$.

An important modelling issue on bounded domains is the choice of the
fractional Laplacian. Several non-equivalent realisations are available and
they generally lead to different boundary-value problems. In particular,
the integral, or restricted, fractional Laplacian involves interactions with
the exterior of the spatial domain, whereas the spectral fractional
Laplacian is defined through the eigenpairs of the classical Dirichlet
Laplacian; see, for example, \cite{Lischke20}. In the present work we employ
the \emph{spectral Dirichlet fractional Laplacian}
\[
(-\Delta_D)^s=\mathcal L^s,
\qquad
\mathcal L=-\frac{d^2}{dx^2},
\]
defined by the fractional powers of the positive Dirichlet Laplacian. This
choice is consistent with the classical homogeneous Dirichlet realisation
and admits a natural semigroup interpretation through subordination of the
Dirichlet heat semigroup \cite{Lischke20}. It is also particularly convenient
for numerical approximation because fractional powers of the corresponding
discrete elliptic operator can be evaluated through its eigenstructure.

For $\Omega=(0,l)$, the Dirichlet eigenfunctions are
$\sin(k\pi x/l)$, and consequently
\[
 (-\Delta_D)^s
 \sin\!\left(\frac{k\pi x}{l}\right)
 =
 \left(\frac{k\pi}{l}\right)^{2s}
 \sin\!\left(\frac{k\pi x}{l}\right).
\]
This property is also useful in the numerical validation carried out below,
since the manufactured solutions can be expressed as finite combinations of
Dirichlet eigenfunctions. The corresponding forcing functions can therefore
be evaluated from the continuous spectral operator independently of the
discrete spatial approximation, avoiding cancellation of spatial
discretisation errors.

The numerical approximation of fractional powers of elliptic operators has
been studied by means of matrix-function techniques, finite-element
approximations, quadrature formulas, and spectral methods; see, for example,
\cite{BonitoPasciak15}. A detailed comparison of different definitions of the
fractional Laplacian on bounded domains is given in \cite{Lischke20}.
Related finite-difference, finite-element, spectral, and quadrature methods
for fractional diffusion equations and other realisations of the fractional
Laplacian can be found in
\cite{DWZ18,DZ19,GSZZ21,HW22,HO14,HS24,WZ22,YCYL25,
Acosta2017,Borthagaray19,CSW21,SWCL24,ZC24}.

Inverse problems with integral overdetermination constitute another important
line of research. Such conditions make it possible to recover unknown
time-dependent quantities from spatially averaged observations rather than
from complete state information. Classical approaches to inverse problems of
this type can be found in \cite{Prilepko2000}, while more recent
reconstruction problems involving time-dependent sources or coefficients are
studied in
\cite{Grim20,VanBock22,Hendy22,Karel22,DurdRah25}.
Inverse problems for fractional diffusion and fractional elliptic equations
have also received considerable attention, including uniqueness and
reconstruction problems, inverse problems with time-dependent coefficients,
and fractional Calder\'on-type problems; see
\cite{Helin20,Ghosh20,Guerngar21,Li23,Kow23,Sahoo23,Fujishiro16,Hazanee13,Hazanee14,Suragan24}.

Particularly relevant is the work of Mamanazarov and Suragan
\cite{Suragan25}, who studied the recovery of a time-dependent source
coefficient in a fractional heat equation with nonlocal data and established
existence and uniqueness of the inverse solution for the regional fractional
Laplacian. The present problem has a similar inverse structure, but we employ
the spectral Dirichlet fractional Laplacian and develop the corresponding
spectral discretisation and numerical reconstruction.

The main contribution of the present work is a complete numerical treatment of this inverse problem for the spectral Dirichlet fractional Laplacian. We derive the continuous Volterra formulation of the inverse problem,
discretise the spatial operator as $A_h=L_h^s$, and combine it with the
Crank--Nicolson method.  Unconditional stability and second-order convergence are proved, and a closed scalar formula for recovering $r(t)$ is derived under a discrete identifiability condition. The effect of noisy integral measurements is also investigated using adaptive regularised differentiation of the measured data $w(t)$.

The remainder of the paper is organised as follows. Section~\ref{2}
introduces the spectral fractional Dirichlet Laplacian and the associated
Hilbert scale.Section~\ref{3} presents the continuous formulation of the inverse problem.  Section~\ref{4}
develops the discrete spectral operator, the Crank--Nicolson scheme, its
stability and convergence analysis, and the coefficient reconstruction
algorithm. Section~\ref{5} presents numerical experiments illustrating the
accuracy of the method, its convergence behaviour, and its sensitivity to
noisy measurements.

\section{Preliminaries}\label{2}
Let
\[
 \Lcal:=-\frac{d^2}{dx^2},\qquad
 D(\Lcal)=H^2(0,l)\cap H_0^1(0,l),
\]
be the positive Dirichlet Laplacian on $L^2(0,l)$. Its eigenpairs are
\begin{equation}\label{eq:eigenpairs}
 \phi_k(x)=\sqrt{\frac{2}{l}}\sin\!\left(\frac{k\pi x}{l}\right),
 \qquad
 \lambda_k=\left(\frac{k\pi}{l}\right)^2,
 \qquad k\ge1,
\end{equation}
and $\{\phi_k\}_{k\ge1}$ is an orthonormal basis of $L^2(0,l)$.

\begin{definition}[Spectral fractional Dirichlet Laplacian]\label{def:spectral-lap}
For $0<s<1$ and
\[
 v=\sum_{k=1}^{\infty}v_k\phi_k,
 \qquad v_k=(v,\phi_k)_{L^2},
\]
define
\begin{equation}\label{eq:spectral-def}
 ( -\Delta_D)^s v=\Lcal^s v
 :=\sum_{k=1}^{\infty}\lambda_k^s v_k\phi_k,
\end{equation}
with domain
\begin{equation}\label{eq:domain-spectral}
 D(\Lcal^s)
 =\left\{v\in L^2(0,l):
 \sum_{k=1}^{\infty}\lambda_k^{2s}|v_k|^2<\infty\right\}.
\end{equation}
\end{definition}

The operator $\Lcal^s$ is positive and self-adjoint. In particular,
\begin{equation}\label{eq:spectral-energy}
 (\Lcal^s v,v)_{L^2}=\|\Lcal^{s/2}v\|_{L^2}^2,
 \qquad v\in D(\Lcal^s).
\end{equation}
For $\sigma\ge0$, it is convenient to introduce the spectral Hilbert scale
\begin{equation}\label{eq:HsD}
 \mathbb H_D^{\sigma}
 :=\left\{v\in L^2(0,l):
 \sum_{k=1}^{\infty}\lambda_k^{\sigma}|v_k|^2<\infty\right\},
 \qquad
 \|v\|_{\mathbb H_D^{\sigma}}^2
 :=\sum_{k=1}^{\infty}\lambda_k^{\sigma}|v_k|^2.
\end{equation}
Thus $D(\Lcal^s)=\mathbb H_D^{2s}$ and
$\|\Lcal^{s/2}v\|_{L^2}=\|v\|_{\mathbb H_D^s}$.

Let $E_s(t)=e^{-t\Lcal^s}$ denote the analytic contraction semigroup generated by $-\Lcal^s$. Its eigenfunction representation is
\begin{equation}\label{eq:semigroup}
 E_s(t)v=\sum_{k=1}^{\infty}e^{-\lambda_k^s t}v_k\phi_k,
 \qquad \|E_s(t)v\|_{L^2}\le \|v\|_{L^2}.
\end{equation}

\begin{definition}[Weak solution]
For a prescribed $r$, a function
\[
 u\in C([0,T];L^2(0,l))\cap L^2(0,T;\mathbb H_D^s),\qquad
 u_t\in L^2(0,T;(\mathbb H_D^s)').
\]
is a weak solution of \eqref{1.1}--\eqref{1.3} if $u(0)=\varphi$ and
\begin{equation}\label{eq:weak}
 \langle u_t,v\rangle
 +(\Lcal^{s/2}u,\Lcal^{s/2}v)_{L^2}
 =r(t)(f(t),v)_{L^2}
\end{equation}
for all admissible test functions $v\in\mathbb H_D^s$, in the usual weak-in-time sense.
\end{definition}

%########################################################################
\section{Continuous formulation of the inverse problem}\label{3}

Before introducing the numerical method, we briefly describe the continuous
reconstruction underlying the inverse problem. Since $-\mathcal L^s$
generates the analytic contraction semigroup
\[
E_s(t)=e^{-t\mathcal L^s},
\]
for a prescribed $r\in L^\infty(0,T)$ the solution of
\eqref{1.1}--\eqref{1.3} is represented by
\begin{equation}\label{eq:mild-solution}
u(t)
=
E_s(t)\varphi
+
\int_0^t
E_s(t-\tau)r(\tau)f(\tau)\,d\tau .
\end{equation}
The corresponding direct problem is well posed by the standard semigroup
theory for positive self-adjoint operators.

Substituting \eqref{eq:mild-solution} into the integral condition
\eqref{1.4} gives
\begin{equation}\label{eq:measurement-volterra}
w(t)
=
(E_s(t)\varphi,\omega)_{L^2}
+
\int_0^t
r(\tau)
(E_s(t-\tau)f(\tau),\omega)_{L^2}
\,d\tau .
\end{equation}

Assume that
\[
\varphi\in D(\mathcal L^s),
\qquad
f\in C([0,T];D(\mathcal L^s)),
\qquad
\omega\in L^2(0,l),
\qquad
w\in C^1[0,T],
\]
together with the compatibility condition
\begin{equation}\label{eq:compatibility}
w(0)=(\varphi,\omega)_{L^2},
\end{equation}
and suppose that
\begin{equation}\label{eq:continuous-nondegeneracy}
|(f(t),\omega)_{L^2}|
\ge\kappa_0>0,
\qquad 0\le t\le T.
\end{equation}

Differentiating \eqref{eq:measurement-volterra} yields
\[
w'(t)
=
-(\mathcal L^sE_s(t)\varphi,\omega)_{L^2}
+
r(t)(f(t),\omega)_{L^2}
-
\int_0^t
r(\tau)
(\mathcal L^sE_s(t-\tau)f(\tau),\omega)_{L^2}
\,d\tau .
\]
Hence the unknown coefficient satisfies the Volterra equation
\begin{equation}\label{eq:continuous-volterra-r}
r(t)
=
g(t)
+
\int_0^t K(t,\tau)r(\tau)\,d\tau ,
\end{equation}
where
\[
g(t)
=
\frac{
w'(t)
+
(\mathcal L^sE_s(t)\varphi,\omega)_{L^2}
}{
(f(t),\omega)_{L^2}
},
\]
and
\[
K(t,\tau)
=
\frac{
(\mathcal L^sE_s(t-\tau)f(\tau),\omega)_{L^2}
}{
(f(t),\omega)_{L^2}
}.
\]
Under \eqref{eq:continuous-nondegeneracy}, standard Volterra theory yields
a unique continuous solution $r$. Together with the compatibility condition,
the corresponding mild solution satisfies the measurement condition
\eqref{1.4}. This continuous formulation provides the basis for the discrete
reconstruction developed below.

\section{Numerical solution of the direct and inverse problems}\label{4}
\subsection{Discrete spectral fractional Laplacian and Crank--Nicolson scheme}
Let $h=l/N$. Set $x_i=ih$ for $i=0,\ldots,N$. For an interior grid vector
$\mathbf v=(v_1,\ldots,v_{N-1})^T$, define
\begin{equation}\label{eq:disc-ip}
 \ip{\mathbf v}{\mathbf z}:=h\sum_{i=1}^{N-1}v_i z_i,
 \qquad \normh{\mathbf v}^2=\ip{\mathbf v}{\mathbf v}.
\end{equation}
The standard second-order Dirichlet difference Laplacian is
\begin{equation}\label{eq:Lh}
 L_h=\frac1{h^2}
 \begin{pmatrix}
 2&-1&&\\
 -1&2&-1&\\
 &\ddots&\ddots&\ddots\\
 &&-1&2
 \end{pmatrix}\in\mathbb R^{(N-1)\times(N-1)}.
\end{equation}
It is diagonalised by the discrete sine basis. More precisely,
\begin{equation}\label{eq:Q}
 Q_{ik}=\sqrt{\frac2N}\sin\!\left(\frac{ik\pi}{N}\right),
 \qquad 1\le i,k\le N-1,
\end{equation}
$Q^TQ=I$, and
\begin{equation}\label{eq:eig-Lh}
 L_h=Q\Lambda_hQ^T,
 \qquad
 \Lambda_h=\operatorname{diag}(\lambda_{1,h},\ldots,\lambda_{N-1,h}),
\end{equation}
with
\begin{equation}\label{eq:lambda-h}
 \lambda_{k,h}=\frac4{h^2}\sin^2\!\left(\frac{k\pi}{2N}\right).
\end{equation}

\begin{definition}[Discrete spectral fractional Laplacian]\label{def:Ah}
We define
\begin{equation}\label{eq:Ah}
 \Ah:=L_h^s=Q\Lambda_h^sQ^T,
 \qquad
 \Lambda_h^s=\operatorname{diag}(\lambda_{1,h}^s,\ldots,\lambda_{N-1,h}^s).
\end{equation}
\end{definition}
The matrix $\Ah$ is symmetric positive definite and, for noninteger $s$, generally dense. Its action need not be formed explicitly: a discrete sine transform, diagonal multiplication by $\lambda_{k,h}^s$, and an inverse sine transform evaluate $\Ah\mathbf v$ efficiently.

Let $t_n=n\tau$, $\tau=T/M$, and let
$\mathbf U^n\approx (u(t_n,x_i))_{i=1}^{N-1}$. Put
\[
 \mathbf F^{n+\frac12}=(f(t_{n+\frac12},x_i))_{i=1}^{N-1},
 \qquad r^{n+\frac12}=r(t_{n+\frac12}).
\]
The Crank--Nicolson scheme for the direct problem is
\begin{equation}\label{fds_full}
 \frac{\mathbf U^{n+1}-\mathbf U^n}{\tau}
 +\Ah\frac{\mathbf U^{n+1}+\mathbf U^n}{2}
 =r^{n+\frac12}\mathbf F^{n+\frac12},
\end{equation}
with
\begin{equation}\label{fds_ic}
 U_i^0=\varphi(x_i),\qquad i=1,\ldots,N-1,
\end{equation}
and $U_0^n=U_N^n=0$. Equivalently,
\begin{equation}\label{eq:CN-matrix}
 \left(I+\frac\tau2\Ah\right)\mathbf U^{n+1}
 =\left(I-\frac\tau2\Ah\right)\mathbf U^n
 +\tau r^{n+\frac12}\mathbf F^{n+\frac12}.
\end{equation}
Define
\begin{equation}\label{eq:A-norm}
 \normAh{\mathbf v}^2:=\ip{\Ah\mathbf v}{\mathbf v},
 \qquad
 \|\mathbf g\|_{A_h^{-1},h}^2:=\ip{\Ah^{-1}\mathbf g}{\mathbf g}.
\end{equation}

\subsection{Stability of the Crank--Nicolson scheme}
\begin{theorem}[Unconditional stability]\label{thm:stability}
Let $\mathbf U^{n+1/2}=(\mathbf U^{n+1}+\mathbf U^n)/2$. Then the following estimates hold for every $\tau>0$.

\emph{(i) Homogeneous case.} If $\mathbf F^{n+1/2}=0$, then
\begin{equation}\label{eq:CN-homogeneous-identity}
 \normh{\mathbf U^{n+1}}^2
 +2\tau\normAh{\mathbf U^{n+1/2}}^2
 =\normh{\mathbf U^n}^2.
\end{equation}

\emph{(ii) $L^2_h$ stability with forcing.}
\begin{equation}\label{eq:L2-growth}
 \normh{\mathbf U^n}
 \le \normh{\mathbf U^0}
 +\sum_{k=0}^{n-1}\tau |r^{k+1/2}|\normh{\mathbf F^{k+1/2}}.
\end{equation}

\emph{(iii) Energy estimate.}
\begin{equation}\label{eq:energy-bound}
 \normh{\mathbf U^n}^2
 +\sum_{k=0}^{n-1}\tau\normAh{\mathbf U^{k+1/2}}^2
 \le \normh{\mathbf U^0}^2
 +\sum_{k=0}^{n-1}\tau |r^{k+1/2}|^2
 \|\mathbf F^{k+1/2}\|_{A_h^{-1},h}^2.
\end{equation}
\end{theorem}
\begin{proof}
Take the discrete inner product of \eqref{fds_full} with
$\mathbf U^{n+1}+\mathbf U^n$. Symmetry of $\Ah$ gives
\begin{equation}\label{eq:energy-core-new}
 \frac{\normh{\mathbf U^{n+1}}^2-\normh{\mathbf U^n}^2}{\tau}
 +2\normAh{\mathbf U^{n+1/2}}^2
 =2r^{n+1/2}\ip{\mathbf F^{n+1/2}}{\mathbf U^{n+1/2}}.
\end{equation}
Part (i) follows immediately. For (ii), one may equivalently diagonalise $\Ah$. On every eigenmode the homogeneous amplification factor is
\[
 g(\mu)=\frac{1-\mu/2}{1+\mu/2},\qquad \mu=\tau\lambda_{k,h}^s\ge0,
\]
so $|g(\mu)|\le1$. The variation-of-constants representation of the discrete system then yields \eqref{eq:L2-growth}. For (iii), use
\[
 |\ip{\mathbf F}{\mathbf V}|
 \le \|\mathbf F\|_{A_h^{-1},h}\normAh{\mathbf V}
\]
and Young's inequality in \eqref{eq:energy-core-new}, followed by summation in $n$.
\end{proof}

\begin{remark}
The stability proof uses only the symmetry and positivity of $\Ah$; hence the Crank--Nicolson part of the analysis is independent of the particular diagonalisation used to apply the spectral fractional operator.
\end{remark}

\subsection{Spatial consistency and convergence}
For a continuous function $v$, write
$\mathbf v_h=(v(x_1),\ldots,v(x_{N-1}))^T$.

\begin{assumption}[Second-order spatial consistency]\label{ass:space}
For all functions appearing in the truncation-error analysis, assume
\begin{equation}\label{eq:space-consistency-ass}
 \normh{\Ah\mathbf v_h-(\Lcal^s v)_h}\le C_v h^2,
\end{equation}
with a constant independent of $h$.
\end{assumption}

The assumption is automatic for the finite sine expansions used in the numerical examples.

\begin{lemma}[Consistency for finite sine expansions]\label{lem:sine-consistency}
Let
\begin{equation}\label{eq:finite-sine}
 v(x)=\sum_{k=1}^{K}a_k\sin\!\left(\frac{k\pi x}{l}\right),
\end{equation}
where $K$ is fixed independently of $h$. Then
\begin{equation}\label{eq:sine-consistency}
 \normh{\Ah\mathbf v_h-(\Lcal^s v)_h}\le C h^2,
\end{equation}
where $C$ depends on $s,l,K$ and the coefficients $a_k$, but not on $h$.
\end{lemma}
\begin{proof}
For each fixed $k\le K$, the nodal sine vector is an eigenvector of $L_h$, hence also of $\Ah$, with eigenvalue $\lambda_{k,h}^s$. The continuous eigenvalue is
$\lambda_k^s=(k\pi/l)^{2s}$. Since
\[
 \lambda_{k,h}
 =\frac4{h^2}\sin^2\!\left(\frac{k\pi h}{2l}\right)
 =\lambda_k-\frac{\lambda_k^2}{12}h^2+O(h^4),
\]
Taylor expansion of $x\mapsto x^s$ at the positive number $\lambda_k$ gives
\[
 \lambda_{k,h}^s=\lambda_k^s+O(h^2).
\]
Summing the finitely many modes proves \eqref{eq:sine-consistency}.
\end{proof}

\begin{remark}
For general smooth data, approximation of fractional powers of elliptic operators can also be analysed through operator-functional calculus; see \cite{BonitoPasciak15}. The finite-mode lemma is stated explicitly because it gives a completely transparent justification for the manufactured solutions used in Section~\ref{5} without importing assumptions from a different spatial discretisation.
\end{remark}

 %------------------------------------------------------------
\begin{theorem}[Second-order convergence of the direct scheme]
\label{thm:direct-convergence}
Assume that the exact solution \(u\) is sufficiently regular and that
\begin{equation}
\sup_{0\le t\le T}
\left(
\|(u_{ttt}(t))_h\|_h
+
\|((rf)_{tt}(t))_h\|_h
\right)
\le C_t,
\label{eq:time-regularity-nodal}
\end{equation}
uniformly with respect to \(h\).
Assume further that the spatial consistency estimate
\begin{equation}
\|A_h v_h-(\mathcal L^s v)_h\|_h
\le C_v h^2
\label{eq:spatial-consistency-thm}
\end{equation}
holds uniformly for
\[
v=u(t,\cdot),\qquad 0\le t\le T,
\]
with a constant independent of \(h\).

Let \(U^n\) be the solution of the Crank--Nicolson scheme
\eqref{fds_full}
with exact initial data
\[
U^0=u_h^0.
\]
Then there exists a constant \(C>0\), independent of \(h\) and \(\tau\),
such that
\begin{equation}
\max_{0\le n\le M}
\|u_h^n-U^n\|_h
\le
C(\tau^2+h^2),
\label{eq:L2-direct-convergence}
\end{equation}
and
\begin{equation}
\max_{0\le n\le M}
\|u_h^n-U^n\|_{A_h,h}
\le
C(\tau^2+h^2).
\label{eq:energy-direct-convergence}
\end{equation}
\end{theorem}

\begin{proof}
Set
\[
e^n=u_h^n-U^n.
\]
We first derive the local consistency relation satisfied by the exact
nodal solution.

For \(t_{n+\frac12}=(t_n+t_{n+1})/2\), the trapezoidal expansion gives
\begin{equation}
\frac{u_h^{n+1}-u_h^n}{\tau}
-
\frac{(u_t)_h^{n+1}+(u_t)_h^n}{2}
=
\theta_1^{n+\frac12},
\qquad
\|\theta_1^{n+\frac12}\|_h
\le C\tau^2,
\label{eq:time-defect-1}
\end{equation}
where \eqref{eq:time-regularity-nodal} has been used.

Similarly,
\begin{equation}
\frac{(rf)_h^{n+1}+(rf)_h^n}{2}
-
r^{n+\frac12}F^{n+\frac12}
=
\theta_2^{n+\frac12},
\qquad
\|\theta_2^{n+\frac12}\|_h
\le C\tau^2.
\label{eq:time-defect-2}
\end{equation}

Averaging the continuous equation
\[
u_t+ \mathcal L^s u=rf
\]
at \(t_n\) and \(t_{n+1}\), and using the spatial consistency
assumption \eqref{eq:spatial-consistency-thm}, we obtain
\begin{equation}
\frac{u_h^{n+1}-u_h^n}{\tau}
+
A_h\frac{u_h^{n+1}+u_h^n}{2}
=
r^{n+\frac12}F^{n+\frac12}
+
\xi^{n+\frac12},
\label{eq:exact-CN-relation-revised}
\end{equation}
where
\begin{equation}
\|\xi^{n+\frac12}\|_h
\le C(\tau^2+h^2).
\label{eq:CN-defect-revised}
\end{equation}

Subtracting the numerical scheme from
\eqref{eq:exact-CN-relation-revised} gives
\begin{equation}
\frac{e^{n+1}-e^n}{\tau}
+
A_h e^{n+\frac12}
=
\xi^{n+\frac12},
\qquad
e^{n+\frac12}:=\frac{e^{n+1}+e^n}{2}.
\label{eq:error-equation-direct}
\end{equation}

Since the Crank--Nicolson propagator
\[
G_h=
\left(I+\frac{\tau}{2}A_h\right)^{-1}
\left(I-\frac{\tau}{2}A_h\right)
\]
satisfies
\[
\|G_h\|_h\le 1,
\]
and
\[
\left\|
\left(I+\frac{\tau}{2}A_h\right)^{-1}
\right\|_h
\le 1,
\]
the discrete variation-of-constants formula yields
\[
\|e^n\|_h
\le
\|e^0\|_h
+
\sum_{k=0}^{n-1}
\tau
\|\xi^{k+\frac12}\|_h.
\]
Because \(e^0=0\), estimate \eqref{eq:CN-defect-revised} gives
\[
\|e^n\|_h
\le
CT(\tau^2+h^2),
\]
which proves \eqref{eq:L2-direct-convergence}.

To obtain the energy estimate, take the discrete inner product of
\eqref{eq:error-equation-direct} with
\(2A_h e^{n+\frac12}\). Since \(A_h\) is symmetric positive definite,
\[
\frac{
\|e^{n+1}\|_{A_h,h}^2
-
\|e^n\|_{A_h,h}^2
}{\tau}
+
2\|A_h e^{n+\frac12}\|_h^2
=
2\langle
\xi^{n+\frac12},
A_h e^{n+\frac12}
\rangle_h.
\]
By Young's inequality,
\[
2\langle
\xi^{n+\frac12},
A_h e^{n+\frac12}
\rangle_h
\le
\|\xi^{n+\frac12}\|_h^2
+
\|A_h e^{n+\frac12}\|_h^2.
\]
Hence
\[
\frac{
\|e^{n+1}\|_{A_h,h}^2
-
\|e^n\|_{A_h,h}^2
}{\tau}
+
\|A_h e^{n+\frac12}\|_h^2
\le
\|\xi^{n+\frac12}\|_h^2.
\]
Summing from \(n=0\) to \(m-1\) and using \(e^0=0\) gives
\[
\|e^m\|_{A_h,h}^2
\le
\tau
\sum_{n=0}^{m-1}
\|\xi^{n+\frac12}\|_h^2
\le
CT(\tau^2+h^2)^2.
\]
Taking square roots proves
\eqref{eq:energy-direct-convergence}.
\end{proof}

\subsection{Numerical algorithm for the inverse problem}
Set
\begin{equation}\label{eq:LR}
 L:=I+\frac\tau2\Ah,\qquad R:=I-\frac\tau2\Ah.
\end{equation}
For a given $\mathbf U^n$, define
\begin{equation}\label{eq:Y-S}
 \mathbf Y=L^{-1}R\mathbf U^n,
 \qquad
 \mathbf S=L^{-1}\mathbf F^{n+1/2},
\end{equation}
and
\begin{equation}\label{eq:V}
 \mathbf V=\frac{\mathbf U^n+\mathbf Y}{2}.
\end{equation}
Then \eqref{eq:CN-matrix} implies
\begin{equation}\label{eq:U-next}
 \mathbf U^{n+1}=\mathbf Y+\tau r^{n+1/2}\mathbf S,
\end{equation}
and hence
\[
 \frac{\mathbf U^{n+1}+\mathbf U^n}{2}
 =\mathbf V+\frac\tau2r^{n+1/2}\mathbf S.
\]

Let $\boldsymbol\omega=(\omega(x_i))_{i=1}^{N-1}$. If an approximation
$z^{n+1/2}\approx w'(t_{n+1/2})$ is available, the continuous measurement identity
\[
 w'(t)+(\Lcal^s u(t),\omega)=r(t)(f(t),\omega)
\]
suggests the discrete midpoint relation
\begin{equation}\label{eq:disc-measure-deriv}
 z^{n+1/2}
 +\ip{\Ah\left(\mathbf V+\frac\tau2r^{n+1/2}\mathbf S\right)}{\boldsymbol\omega}
 =r^{n+1/2}\ip{\mathbf F^{n+1/2}}{\boldsymbol\omega}.
\end{equation}
Since $L\mathbf S=\mathbf F^{n+1/2}$,
\begin{equation}\label{eq:den-simplify}
 \mathbf F^{n+1/2}=\mathbf S+\frac\tau2\Ah\mathbf S,
\end{equation}
so the coefficient multiplying $r^{n+1/2}$ in \eqref{eq:disc-measure-deriv} reduces exactly to
\begin{equation}\label{eq:d-h}
 d_h^{n+1/2}:=\ip{\mathbf S}{\boldsymbol\omega}.
\end{equation}
Therefore
\begin{equation}\label{eq:r-closed}
 r_h^{n+1/2}
 =\frac{z^{n+1/2}+\ip{\Ah\mathbf V}{\boldsymbol\omega}}
 {\ip{\mathbf S}{\boldsymbol\omega}}
\end{equation}

If only sampled values $w^n=w(t_n)$ are available, one may take
\begin{equation}\label{eq:z-difference}
 z^{n+1/2}=D_\tau w^{n+1/2}:=\frac{w^{n+1}-w^n}{\tau},
\end{equation}
which is a second-order approximation to $w'(t_{n+1/2})$ for smooth exact data. For noisy data, however, direct use of \eqref{eq:z-difference} should be replaced by a regularised differentiation or smoothing procedure; see Section~\ref{subsec:noisy-data}.

%------------------------------------------------------------
\begin{lemma}[Discrete identifiability]
\label{lem:discrete-identifiability}
Assume the continuous nondegeneracy condition
\begin{equation}
|(f(t),\omega)_{L^2}|
\ge
\kappa_0>0,
\qquad
0\le t\le T.
\label{eq:continuous-nondegeneracy-recall}
\end{equation}
Suppose further that
\begin{equation}
\left|
\langle F^{n+\frac12},\omega\rangle_h
-
(f(t_{n+\frac12}),\omega)_{L^2}
\right|
\le
C_Qh^2,
\label{eq:quadrature-fomega}
\end{equation}
uniformly in \(n\), and that
\begin{equation}
\sup_{n,h}\|F^{n+\frac12}\|_h\le C_F,
\qquad
\sup_h\|A_h\omega\|_h\le C_\omega.
\label{eq:identifiability-bounds}
\end{equation}

Let
\[
S^{n+\frac12}
=
\left(I+\frac{\tau}{2}A_h\right)^{-1}
F^{n+\frac12},
\]
and define
\[
d_h^{n+\frac12}
=
\langle S^{n+\frac12},\omega\rangle_h.
\]
Then
\begin{equation}
\left|
d_h^{n+\frac12}
-
(f(t_{n+\frac12}),\omega)_{L^2}
\right|
\le
C_Qh^2
+
\frac{\tau}{2}C_FC_\omega.
\label{eq:discrete-denominator-consistency}
\end{equation}
Consequently, if \(h\) and \(\tau\) are sufficiently small so that
\begin{equation}
C_Qh^2
+
\frac{\tau}{2}C_FC_\omega
\le
\frac{\kappa_0}{2},
\label{eq:mesh-identifiability-condition}
\end{equation}
then
\begin{equation}
|d_h^{n+\frac12}|
\ge
\frac{\kappa_0}{2}
=:\kappa
>0,
\qquad
n=0,\ldots,M-1.
\label{eq:discrete-identifiability}
\end{equation}
\end{lemma}

\begin{proof}
By definition,
\[
\left(I+\frac{\tau}{2}A_h\right)
S^{n+\frac12}
=
F^{n+\frac12}.
\]
Therefore
\[
F^{n+\frac12}
=
S^{n+\frac12}
+
\frac{\tau}{2}A_hS^{n+\frac12},
\]
and hence
\[
d_h^{n+\frac12}
=
\langle F^{n+\frac12},\omega\rangle_h
-
\frac{\tau}{2}
\langle A_hS^{n+\frac12},\omega\rangle_h.
\]
Since \(A_h\) is symmetric,
\[
\langle A_hS^{n+\frac12},\omega\rangle_h
=
\langle S^{n+\frac12},A_h\omega\rangle_h.
\]
Moreover,
\[
\left\|
\left(I+\frac{\tau}{2}A_h\right)^{-1}
\right\|_h
\le1,
\]
so
\[
\|S^{n+\frac12}\|_h
\le
\|F^{n+\frac12}\|_h
\le C_F.
\]
Consequently,
\[
\left|
d_h^{n+\frac12}
-
\langle F^{n+\frac12},\omega\rangle_h
\right|
\le
\frac{\tau}{2}C_FC_\omega.
\]
Combining this estimate with
\eqref{eq:quadrature-fomega} proves
\eqref{eq:discrete-denominator-consistency}.

Finally,
\[
|d_h^{n+\frac12}|
\ge
|(f(t_{n+\frac12}),\omega)_{L^2}|
-
C_Qh^2
-
\frac{\tau}{2}C_FC_\omega.
\]
Using \eqref{eq:continuous-nondegeneracy-recall} and
\eqref{eq:mesh-identifiability-condition} yields
\[
|d_h^{n+\frac12}|
\ge
\frac{\kappa_0}{2}.
\]
\end{proof}

\begin{remark}
Lemma~\ref{lem:discrete-identifiability} shows that the discrete
identifiability condition is not an independent assumption.
It follows from the continuous nondegeneracy condition
\eqref{eq:continuous-nondegeneracy-recall} whenever the spatial and
temporal meshes are sufficiently fine. In practical computations,
the quantity
\[
|d_h^{n+\frac12}|
=
|\langle S^{n+\frac12},\omega\rangle_h|
\]
should nevertheless be monitored at every time step.
\end{remark}

 The preceding relations lead to a sequential reconstruction procedure
for the unknown coefficient and the state. At each time level, the
known state $\mathbf U^n$ is first propagated through the homogeneous
Crank--Nicolson step, while the response $ S$ associated with
the source term is computed from the same system matrix. The
coefficient $r_h^{n+\frac12}$ is then recovered from the scalar
formula \eqref{eq:r-closed}, after which the state is updated to
$\mathbf U^{n+1}$. The complete reconstruction procedure is summarised
in Algorithm~\ref{alg:CN-SFL}.

\begin{algorithm}[H]
\caption{CN--SFL algorithm for simultaneous recovery of $r(t)$ and $u(t,x)$}
\label{alg:CN-SFL}
\begin{algorithmic}[1]
\Require Grids $x_i=ih$, $t_n=n\tau$; data $\varphi$, $f$, $\omega$; measured values $w^n$.
\State Construct the eigenvalues of $L_h$ and the spectral
representation of $A_h=L_h^s$ according to \eqref{eq:Ah};
apply the resulting operators by the discrete sine transform.
\State Set $\mathbf U^0=(\varphi(x_i))_{i=1}^{N-1}$ and $ \omega=(\omega(x_i))_{i=1}^{N-1}$.
\State Set $L=I+\frac\tau2\Ah$ and $R=I-\frac\tau2\Ah$.
\State Construct $z^{n+1/2}$ from the measurements $w^n$; for noisy data use the regularised differentiation procedure described in Section~\ref{subsec:noisy-data}.
\For{$n=0$ to $M-1$}
 \State Assemble $  F^{n+1/2}$.
 \State Solve $L\mathbf Y=R\mathbf U^n$ and $L  S=  F^{n+1/2}$.
 \State Set $\mathbf V=(\mathbf U^n+\mathbf Y)/2$.
 \State Compute $d_h^{n+\frac12}
=\langle\mathbf S,\boldsymbol\omega\rangle_h$.
\State Monitor $|d_h^{n+\frac12}|$; by
Lemma~\ref{lem:discrete-identifiability}, this quantity is uniformly
separated from zero for sufficiently fine meshes.
 \State Recover $r_h^{n+1/2}$ from \eqref{eq:r-closed}.
 \State Update $\mathbf U^{n+1}=\mathbf Y+\tau r_h^{n+1/2}  S$.
\EndFor
\State \Return $\{\mathbf U^n\}_{n=0}^M$ and $\{r_h^{n+1/2}\}_{n=0}^{M-1}$.
\end{algorithmic}
\end{algorithm}

Thus, once the derivative approximation $z^{n+\frac12}$ is available,
the inverse problem is reduced at each time step to two linear solves
with the same matrix $L$ and one scalar coefficient evaluation.

\subsection{Conditional convergence of the inverse reconstruction}
The following theorem makes explicit the additional assumptions needed to transfer direct-scheme convergence to the recovered coefficient.

%------------------------------------------------------------
\begin{theorem}[Conditional convergence of the coupled reconstruction]
\label{thm:inverse-convergence}
Assume the hypotheses of
Theorem~\ref{thm:direct-convergence} and
Lemma~\ref{lem:discrete-identifiability}.
Assume in addition that
\begin{equation}
\sup_h\|A_h\omega\|_h\le C_\omega,
\qquad
\sup_{n,h}\|F^{n+\frac12}\|_h\le C_F,
\label{eq:inverse-uniform-bounds}
\end{equation}
and that
\[
w\in C^3[0,T].
\]

Suppose that the derivative approximation constructed from the
measurement data satisfies
\begin{equation}
\left|
z^{n+\frac12}
-
w'(t_{n+\frac12})
\right|
\le
\eta,
\qquad
n=0,\ldots,M-1.
\label{eq:derivative-error-revised}
\end{equation}
Assume also that the nodal quadrature of the differentiated
measurement satisfies
\begin{equation}
\left|
w'(t)
-
\left\langle
(u_t(t,x_i))_{i=1}^{N-1},
\omega
\right\rangle_h
\right|
\le
C_Qh^2,
\qquad
0\le t\le T.
\label{eq:measurement-quadrature-revised}
\end{equation}

Then, for sufficiently small \(h\) and \(\tau\), there exists a
constant \(C>0\), independent of \(h\), \(\tau\), and \(\eta\), such
that
\begin{equation}
\max_{0\le n\le M}
\|u_h^n-U^n\|_h
+
\max_{0\le n\le M-1}
\left|
r(t_{n+\frac12})
-
r_h^{n+\frac12}
\right|
\le
C(\tau^2+h^2+\eta).
\label{eq:coupled-convergence-revised}
\end{equation}
\end{theorem}

\begin{proof}
Set
\[
e^n=u_h^n-U^n,
\qquad
\rho^{n+\frac12}
=
r(t_{n+\frac12})-r_h^{n+\frac12}.
\]

For brevity, write
\[
r_*^{n+\frac12}:=r(t_{n+\frac12}),
\qquad
L=I+\frac{\tau}{2}A_h,
\qquad
R=I-\frac{\tau}{2}A_h.
\]
The exact nodal solution satisfies
\begin{equation}
Lu_h^{n+1}
=
Ru_h^n
+
\tau r_*^{n+\frac12}F^{n+\frac12}
+
\tau\xi^{n+\frac12},
\label{eq:exact-state-LR}
\end{equation}
where
\begin{equation}
\|\xi^{n+\frac12}\|_h
\le
C(\tau^2+h^2)
\label{eq:inverse-defect-bound}
\end{equation}
by Theorem~\ref{thm:direct-convergence}.

Define
\[
Y_*^n=L^{-1}Ru_h^n,
\qquad
S^{n+\frac12}=L^{-1}F^{n+\frac12},
\qquad
V_*^n=\frac{u_h^n+Y_*^n}{2}.
\]
Equation \eqref{eq:exact-state-LR} gives
\[
u_h^{n+1}
=
Y_*^n
+
\tau r_*^{n+\frac12}S^{n+\frac12}
+
\tau L^{-1}\xi^{n+\frac12}.
\]
Therefore
\begin{equation}
\frac{u_h^{n+1}+u_h^n}{2}
=
V_*^n
+
\frac{\tau}{2}
r_*^{n+\frac12}S^{n+\frac12}
+
\frac{\tau}{2}L^{-1}\xi^{n+\frac12}.
\label{eq:exact-midpoint-state}
\end{equation}

We next derive the corresponding scalar identity for the exact
coefficient. Taking the discrete inner product of
\eqref{eq:exact-CN-relation-revised} with \(\omega\) gives
\begin{equation}
\left\langle
\frac{u_h^{n+1}-u_h^n}{\tau},
\omega
\right\rangle_h
+
\left\langle
A_h\frac{u_h^{n+1}+u_h^n}{2},
\omega
\right\rangle_h
=
r_*^{n+\frac12}
\langle F^{n+\frac12},\omega\rangle_h
+
\langle\xi^{n+\frac12},\omega\rangle_h.
\label{eq:exact-measurement-discrete}
\end{equation}

Since
\[
\frac{1}{\tau}
\left\langle
u_h^{n+1}-u_h^n,\omega
\right\rangle_h
=
\frac{1}{\tau}
\int_{t_n}^{t_{n+1}}
\langle (u_t(t,x_i))_i,\omega\rangle_h\,dt,
\]
assumption \eqref{eq:measurement-quadrature-revised} implies
\[
\left|
\frac{1}{\tau}
\left\langle
u_h^{n+1}-u_h^n,\omega
\right\rangle_h
-
\frac{w(t_{n+1})-w(t_n)}{\tau}
\right|
\le
C h^2.
\]
Since \(w\in C^3[0,T]\),
\[
\frac{w(t_{n+1})-w(t_n)}{\tau}
=
w'(t_{n+\frac12})
+
O(\tau^2).
\]
Hence
\begin{equation}
\left|
\frac{1}{\tau}
\left\langle
u_h^{n+1}-u_h^n,\omega
\right\rangle_h
-
w'(t_{n+\frac12})
\right|
\le
C(\tau^2+h^2).
\label{eq:measurement-difference-bound}
\end{equation}

Using \eqref{eq:exact-midpoint-state},
the symmetry of \(A_h\), and
\[
F^{n+\frac12}
=
S^{n+\frac12}
+
\frac{\tau}{2}A_hS^{n+\frac12},
\]
equation \eqref{eq:exact-measurement-discrete} can be rearranged as
\begin{equation}
r_*^{n+\frac12}
d_h^{n+\frac12}
=
w'(t_{n+\frac12})
+
\langle A_hV_*^n,\omega\rangle_h
+
\zeta^{n+\frac12},
\label{eq:exact-coefficient-identity}
\end{equation}
where
\begin{equation}
|\zeta^{n+\frac12}|
\le
C(\tau^2+h^2).
\label{eq:zeta-bound}
\end{equation}
Indeed, the terms contributing to \(\zeta^{n+\frac12}\) are
the measurement approximation error
\eqref{eq:measurement-difference-bound},
the defect term
\(\langle\xi^{n+\frac12},\omega\rangle_h\),
and
\[
\frac{\tau}{2}
\langle
A_hL^{-1}\xi^{n+\frac12},
\omega
\rangle_h.
\]
For the last term, symmetry gives
\[
\left|
\left\langle
A_hL^{-1}\xi^{n+\frac12},
\omega
\right\rangle_h
\right|
=
\left|
\left\langle
L^{-1}\xi^{n+\frac12},
A_h\omega
\right\rangle_h
\right|
\le
C_\omega
\|\xi^{n+\frac12}\|_h,
\]
because \(\|L^{-1}\|_h\le1\).

The numerical coefficient satisfies
\begin{equation}
r_h^{n+\frac12}
d_h^{n+\frac12}
=
z^{n+\frac12}
+
\langle A_hV^n,\omega\rangle_h,
\label{eq:numerical-coefficient-identity}
\end{equation}
where
\[
V^n
=
\frac{U^n+L^{-1}RU^n}{2}.
\]
Subtracting
\eqref{eq:numerical-coefficient-identity}
from \eqref{eq:exact-coefficient-identity} gives
\[
\rho^{n+\frac12}d_h^{n+\frac12}
=
w'(t_{n+\frac12})-z^{n+\frac12}
+
\langle A_h(V_*^n-V^n),\omega\rangle_h
+
\zeta^{n+\frac12}.
\]
By symmetry,
\[
\left|
\langle A_h(V_*^n-V^n),\omega\rangle_h
\right|
\le
\|V_*^n-V^n\|_h
\|A_h\omega\|_h.
\]
Moreover,
\[
V_*^n-V^n
=
\frac12
\left(
I+L^{-1}R
\right)e^n.
\]
Since
\[
\|L^{-1}R\|_h\le1,
\]
we obtain
\[
\|V_*^n-V^n\|_h
\le
\|e^n\|_h.
\]
Lemma~\ref{lem:discrete-identifiability} gives
\[
|d_h^{n+\frac12}|
\ge\kappa>0.
\]
Consequently,
\begin{equation}
|\rho^{n+\frac12}|
\le
C\left(
\|e^n\|_h
+
\tau^2+h^2+\eta
\right).
\label{eq:rho-estimate-revised}
\end{equation}

It remains to estimate the state error. Subtracting the numerical
state equation from \eqref{eq:exact-state-LR} yields
\[
Le^{n+1}
=
Re^n
+
\tau\rho^{n+\frac12}F^{n+\frac12}
+
\tau\xi^{n+\frac12}.
\]
Since
\[
\|L^{-1}R\|_h\le1,
\qquad
\|L^{-1}\|_h\le1,
\]
and \(\|F^{n+\frac12}\|_h\le C_F\),
\[
\|e^{n+1}\|_h
\le
\|e^n\|_h
+
C\tau|\rho^{n+\frac12}|
+
C\tau(\tau^2+h^2).
\]
Using \eqref{eq:rho-estimate-revised},
\[
\|e^{n+1}\|_h
\le
(1+C\tau)\|e^n\|_h
+
C\tau(\tau^2+h^2+\eta).
\]
Since \(e^0=0\), the discrete Gr\"onwall inequality yields
\[
\max_{0\le n\le M}\|e^n\|_h
\le
C(\tau^2+h^2+\eta).
\]
Substitution into \eqref{eq:rho-estimate-revised} gives
\[
\max_{0\le n\le M-1}
|\rho^{n+\frac12}|
\le
C(\tau^2+h^2+\eta).
\]
This proves \eqref{eq:coupled-convergence-revised}.
\end{proof}

%------------------------------------------------------------
\begin{remark}[Raw noise and numerical differentiation]
\label{rem:raw-noise}
Suppose that the available measurements are
\[
w^{\delta,n}
=
w(t_n)+e_n,
\qquad
\varepsilon_\infty
:=
\max_{0\le n\le M}|e_n|.
\]
If the derivative is approximated directly by
\[
D_\tau w^{\delta,n+\frac12}
=
\frac{w^{\delta,n+1}-w^{\delta,n}}{\tau},
\]
then
\[
\left|
D_\tau w^{\delta,n+\frac12}
-
D_\tau w^{n+\frac12}
\right|
=
\frac{|e_{n+1}-e_n|}{\tau}
\le
\frac{2\varepsilon_\infty}{\tau}.
\]
Thus, decreasing the time step does not necessarily improve the
coefficient reconstruction when noisy measurements are differentiated
directly; on the contrary, the noise contribution may grow like
\(\varepsilon_\infty/\tau\).

Theorem~\ref{thm:inverse-convergence} therefore expresses the
reconstruction error in terms of the error
\[
\eta_\delta
=
\max_n
\left|
z^{\delta,n+\frac12}
-
w'(t_{n+\frac12})
\right|
\]
of a regularised derivative approximation, rather than directly in
terms of the raw measurement perturbation.
No convergence rate of \(\eta_\delta\) with respect to the original
noise level \(\delta\) is assumed or asserted here.
\end{remark}

%########################################################################
\section{Numerical Experiments}\label{5}
We consider two manufactured examples for the spectral fractional Laplacian. Since all exact states are finite sums of sine eigenfunctions, the forcing term is evaluated analytically from \eqref{eq:spectral-def}; it is \emph{not} generated with the discrete matrix $\Ah$. This avoids cancellation of the spatial discretisation error.

Unless otherwise stated, $l=T=1$. The eigenfunction identity
\[
 ( -\Delta_D)^s\sin(k\pi x)=(k\pi)^{2s}\sin(k\pi x)
\]
is used throughout.

\noindent\textbf{Example 1.}
Let
\[
 u(t,x)=a(t)\sin(\pi x)+b(t)\sin(3\pi x),\,\,
 a(t)=1+t^2+s\sin t,\,\, b(t)=s t e^{-t},
\]
\[
 r(t)=1+\frac{s}{2}(1+\cos t),\,\, \omega(x)=\sin(\pi x).
\]
By orthogonality,
\[
 w(t)=\int_0^1u(t,x)\omega(x)\,dx=\frac12a(t),
\]
and
\[
 f(t,x)=\frac{1}{r(t)}\left[
 \bigl(a'(t)+\pi^{2s}a(t)\bigr)\sin(\pi x)
 +\bigl(b'(t)+(3\pi)^{2s}b(t)\bigr)\sin(3\pi x)
 \right].
\]
This formula is an exact continuous forcing for the spectral operator.

\noindent\textbf{Example 2.}
We take
\[
 u(t,x)=\cos t\,\sin(\pi x)+s t e^{-t}\sin(2\pi x), \,\,
 r(t)=1+\sin t,\,\, \omega(x)=x(1-x).
\]
The smooth observation weight $\omega(x)=x(1-x)$ is compatible with the regularity assumptions of the analysis. Since $x(1-x)$ is symmetric about $x=1/2$, its integral against $\sin(2\pi x)$ vanishes, while
$\int_0^1x(1-x)\sin(\pi x)\,dx=4/\pi^3$. Hence
$ w(t)=\frac4{\pi^3}\cos t.$

The forcing is
\[
 f(t,x)=\frac1{1+\sin t}\Big[
 \bigl(-\sin t+\pi^{2s}\cos t\bigr)\sin(\pi x)
 +s e^{-t}\bigl(1-t+(2\pi)^{2s}t\bigr)\sin(2\pi x)
 \Big].
\]

\begin{figure}[H]
\centering
\includegraphics[width=0.8\textwidth]{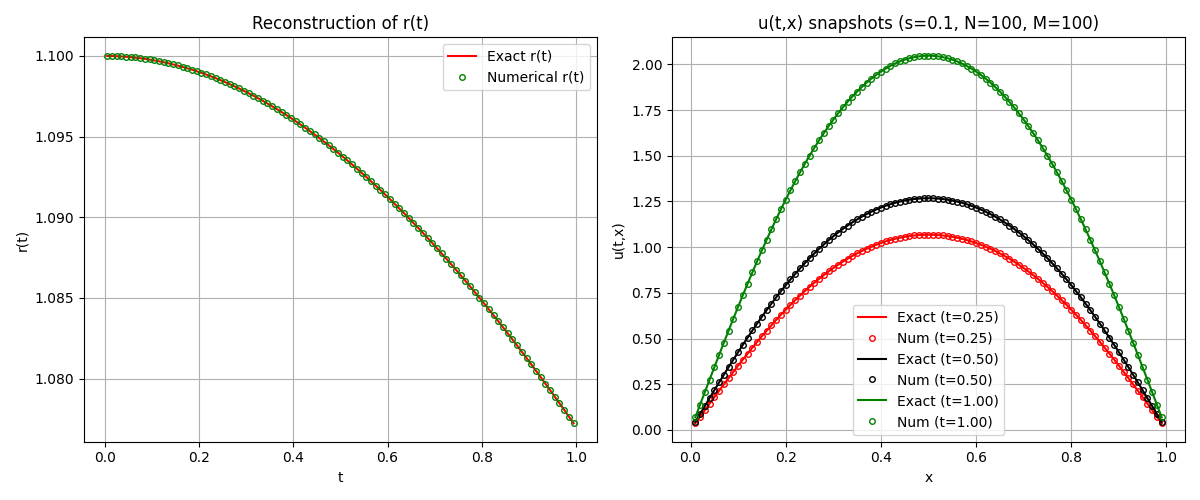}
\caption{Analytical and numerical solutions of $r(t)$ and $u(t,x)$ for Example~1 at $T=1$.}
\label{fig1}
\end{figure}

Figure~\ref{fig1} compares the analytical and reconstructed solutions of
Example~1 for $s=0.1$ and $N=M=100$.

\begin{figure}[H]
\centering
\includegraphics[width=0.8\textwidth]{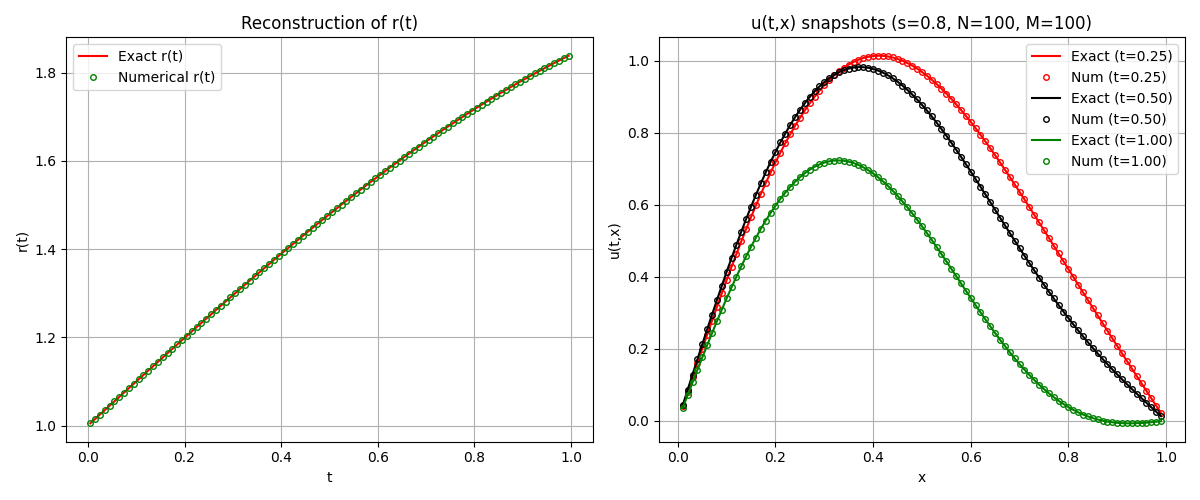}
\caption{Analytical and numerical solutions of $r(t)$ and $u(t,x)$ for Example~2 at $T=1$.}
\label{fig2}
\end{figure}

Figure~\ref{fig2} shows the corresponding comparison for Example~2 with
$s=0.8$ and $N=M=100$.

\subsection{Temporal and spatial convergence tests}

For the convergence tests, the errors in the state variable are evaluated
at the final time $T$. We use the discrete maximum norm
\[
E_{\infty,u}
:=
\max_{1\le i\le N-1}
\left|
U_i^M-u(T,x_i)
\right|,
\]
and the discrete $L^2$ norm
\[
E_{2,u}
:=
\left(
h\sum_{i=1}^{N-1}
\left|
U_i^M-u(T,x_i)
\right|^2
\right)^{1/2}.
\]
For the reconstructed coefficient, the error is measured at the midpoint
time levels by
\[
E_{\infty,r}
:=
\max_{0\le n\le M-1}
\left|
r_h^{n+\frac12}
-r(t_{n+\frac12})
\right|.
\]

In these noise-free convergence experiments, the exact integral data
$w(t_n)$ are used without perturbation. The derivative values required by
the reconstruction formula are approximated by
\[
z^{n+\frac12}
=
D_\tau w^{n+\frac12}
:=
\frac{w(t_{n+1})-w(t_n)}{\tau}.
\]
For sufficiently smooth $w$, this is a second-order approximation at the
midpoint,
\[
D_\tau w^{n+\frac12}
=
w'(t_{n+\frac12})+O(\tau^2).
\]
Thus, no regularisation is required in Tables~\ref{tab:integration-tau}
and~\ref{tab:integration-h}; regularised
differentiation is used only in the noisy-data experiments of
Section~\ref{subsec:noisy-data}.

For successive mesh refinements, the observed convergence order is computed
as
\[
\mathrm{CO}
=
\log_2\left(\frac{E_{\mathrm{coarse}}}
{E_{\mathrm{fine}}}\right).
\]
To assess the temporal order, the spatial mesh is fixed at $h=1/25600$ so that the $O(h^2)$ spatial error is negligible while $\tau$ is refined. We take $s=0.8$. The corresponding temporal errors and observed
convergence orders are reported in Table~\ref{tab:integration-tau}.

\begin{table}[H]
\centering
\caption{Temporal convergence for varying $\tau$ with fixed $h=1/25600$ and $s=0.8$.}
\label{tab:integration-tau}
\begin{tabular}{c|c|c|c|c|c|c}
\hline
$\tau$ & $E_{\infty,u}$  & $\mathrm{CO}$ & $E_{2,u}$ & $\mathrm{CO}$  & $E_{\infty,r}$ & $\mathrm{CO}$ \\
\hline
1/50  &2.176e-05 &--&1.538e-05&--&1.540e-04&--\\
1/100 &5.441e-06&1.999&3.848e-06&1.999&3.876e-05&1.990\\
1/200 &1.362e-06&1.998&9.632e-07&1.998& 9.724e-06&1.995\\
1/400 &3.423e-07&1.992&2.421e-07&1.992&2.434e-06&1.998\\
1/800 &8.739e-08&1.970&6.179e-08&1.970&6.078e-07&2.002\\
\hline
\end{tabular}
\end{table}

To assess the spatial order, the time step is fixed at $\tau=1/12800$, which makes the $O(\tau^2)$ temporal contribution negligible relative to the spatial error over the reported meshes. We take $s=0.1$. The corresponding spatial errors and observed convergence orders are
presented in Table~\ref{tab:integration-h}.
\begin{table}[H]
\centering
\caption{Spatial convergence for varying $h$ with fixed $\tau=1/12800$ and $s=0.1$.}
\label{tab:integration-h}
\begin{tabular}{c|c|c|c|c|c|c}
\hline
$h$ & $E_{\infty,u}$ & $\mathrm{CO}$  & $E_{2,u}$ & $\mathrm{CO}$  & $E_{\infty,r}$ & $\mathrm{CO}$ \\
\hline
1/100 &1.594e-06&--&1.127e-06&--&8.379e-06&--\\
1/200 &3.984e-07&2.000&2.817e-07&2.000&2.093e-06&2.001\\
1/400 &9.965e-08&1.999&7.046e-08&1.999&5.222e-07&2.003\\
1/800 &2.496e-08&1.998&1.765e-08&1.998&1.294e-07&2.013\\
\hline
\end{tabular}
\end{table}
As seen from Tables~\ref{tab:integration-tau}
and~\ref{tab:integration-h}, the observed temporal and spatial
convergence orders are approximately two. In particular, the
second-order behaviour of $E_{2,u}$ and $E_{\infty,r}$ is consistent
with Theorems~\ref{thm:direct-convergence}
and~\ref{thm:inverse-convergence}. The corresponding second-order
behaviour observed for $E_{\infty,u}$ is an additional numerical
observation.

%------------------------------------------------------------
\subsection{Sensitivity to noisy data}
\label{subsec:noisy-data}

We next investigate numerically the sensitivity of the inverse
reconstruction to perturbations in the measured integral data
\[
w(t)
=
\int_0^1 u(t,x)\omega(x)\,dx.
\]
Noise is introduced only into the measured quantity \(w\).
More precisely, at the discrete time levels \(t_n=n\tau\) we set
\begin{equation}
w^\delta(t_n)
=
w(t_n)+e_n,
\qquad
n=0,\ldots,M,
\label{eq:noisy-measurements}
\end{equation}
where
$e=(e_0,\ldots,e_M)^T$
is a zero-mean Gaussian perturbation vector rescaled so that
\begin{equation}
\frac{
\|w^\delta-w\|_{\ell^2}
}{
\|w\|_{\ell^2}
}
=
\delta,
\label{eq:relative-noise-level}
\end{equation}
with 
$\delta\in\{0.01,0.03,0.05\}.$
Thus the experiments correspond to \(1\%\), \(3\%\), and \(5\%\)
relative perturbations of the measured data.

The coefficient reconstruction formula requires an approximation of
$z(t)=w'(t).$

Direct differentiation of the noisy measurements is inappropriate,
because
\[
\frac{
w^\delta(t_{n+1})-w^\delta(t_n)
}{\tau}
-
\frac{
w(t_{n+1})-w(t_n)
}{\tau}
=
\frac{e_{n+1}-e_n}{\tau},
\]
and therefore
\[
\left|
\frac{e_{n+1}-e_n}{\tau}
\right|
\le
\frac{2\varepsilon_\infty}{\tau},
\qquad
\varepsilon_\infty
=
\max_n|e_n|.
\]
Hence differentiation may strongly amplify high-frequency
measurement noise when \(\tau\) is small.

To obtain a stable derivative approximation, the noisy data are first
processed by a Savitzky--Golay local polynomial smoother of degree
\(p=2\). Let $S_q\in\mathbb R^{(M+1)\times(M+1)}$
denote the linear smoothing matrix corresponding to an odd window
length \(q\). We consider the candidate set
$q\in\{7,9,\ldots,81\}.$

For each \(q\), the smoothed data are
$w_q^\delta=S_qw^\delta,$
and the residual is defined by
\begin{equation}
R_q
=
\|(I-S_q)w^\delta\|_{\ell^2}.
\label{eq:SG-residual}
\end{equation}

For an independent Gaussian noise vector with covariance
\(\sigma^2I\), the expected squared residual noise satisfies
\[
\mathbb E
\|(I-S_q)e\|_{\ell^2}^2
=
\sigma^2
\operatorname{tr}
\left(
(I-S_q)^T(I-S_q)
\right).
\]
Accordingly, we introduce the residual degrees-of-freedom quantity
\begin{equation}
\nu_q
:=
\operatorname{tr}
\left(
(I-S_q)^T(I-S_q)
\right).
\label{eq:SG-residual-dof}
\end{equation}
Notice that, in general,
$\nu_q
\neq
M+1-\operatorname{tr}(S_q),$
because a Savitzky--Golay smoothing matrix is not generally an
orthogonal projection.

Since the synthetic noise level \(\delta\) is prescribed relative to
the unknown exact signal, we estimate its discrete
\(\ell^2\)-magnitude from the noisy measurements. Under the
approximate orthogonality relation
$\langle w,e\rangle_{\ell^2}\approx0,$
we have
\[
\|w^\delta\|_{\ell^2}^2
\approx
\|w\|_{\ell^2}^2+\|e\|_{\ell^2}^2
=
(1+\delta^2)\|w\|_{\ell^2}^2.
\]
Hence
\begin{equation}
\varepsilon_{2,\delta}
:=
\frac{\delta}{\sqrt{1+\delta^2}}
\|w^\delta\|_{\ell^2}
\label{eq:estimated-noise-norm}
\end{equation}
is used as an estimate of \(\|e\|_{\ell^2}\).

Motivated by the Gaussian noise model, the discrepancy target
associated with the window \(q\) is defined as
\begin{equation}
R_q^{\mathrm{tar}}
=
\varepsilon_{2,\delta}
\sqrt{
\frac{\nu_q}{M+1}
}.
\label{eq:corrected-discrepancy-target}
\end{equation}
The corresponding relative discrepancy score is
\begin{equation}
D_q
=
\frac{
|R_q-R_q^{\mathrm{tar}}|
}{
R_q^{\mathrm{tar}}
}.
\label{eq:relative-discrepancy-score}
\end{equation}

The smoothing window is selected by minimising \(D_q\) over all
admissible odd window lengths. If several windows have discrepancy
scores within \(5\times10^{-3}\) of the minimum value, the smallest
such window is selected in order to avoid unnecessary oversmoothing.

Once the window has been selected, the regularised derivative
\(z^\delta\) is obtained by differentiating the same local polynomial
fit. In particular, the midpoint values
\[
z^{\delta,n+\frac12}
\approx
w'(t_{n+\frac12})
\]
are inserted into the reconstruction formula
\begin{equation}
r_\delta^{n+\frac12}
=
\frac{
z^{\delta,n+\frac12}
+
\langle A_hV^n,\omega\rangle_h
}{
\langle S^{n+\frac12},\omega\rangle_h
}.
\label{eq:noisy-reconstruction-formula}
\end{equation}

Theorem~\ref{thm:inverse-convergence} shows that the quantity directly
controlling the coefficient error is the derivative reconstruction
error
\begin{equation}
\eta_\delta
:=
\max_{0\le n\le M-1}
\left|
z^{\delta,n+\frac12}
-
w'(t_{n+\frac12})
\right|.
\label{eq:eta-delta}
\end{equation}
Thus, for the manufactured examples considered here, $\eta_\delta$
can be evaluated explicitly because the exact function \(w'(t)\) is
known. This quantity is used only for assessing the numerical
experiment; it is not available in an actual inverse problem.

Provided that the discrete identifiability condition of
Lemma~\ref{lem:discrete-identifiability} remains satisfied, the
conditional estimate of
Theorem~\ref{thm:inverse-convergence} takes the form
\begin{equation}
\max_{0\le n\le M}
\|u_h^n-U_\delta^n\|_h
+
\max_{0\le n\le M-1}
\left|
r(t_{n+\frac12})
-
r_\delta^{n+\frac12}
\right|
\le
C\left(
\tau^2+h^2+\eta_\delta
\right).
\label{eq:noisy-conditional-estimate}
\end{equation}
We emphasise that
\eqref{eq:noisy-conditional-estimate} is a conditional error estimate
with respect to the differentiated-data error \(\eta_\delta\).
The present analysis does not establish a convergence rate
\(\eta_\delta=\eta_\delta(\delta)\) for the Savitzky--Golay
differentiation procedure as \(\delta\to0\). The numerical experiments
below are therefore intended to demonstrate the empirical robustness
of the reconstruction with respect to noisy measurements rather than
to constitute a regularisation convergence theorem.

For the noisy-data computations we take
$N=M=100,
\,
s=0.9.$

The same Gaussian perturbation realisation is used for all three
noise levels, with only its amplitude rescaled according to
\eqref{eq:relative-noise-level}. The selected window length is
determined independently for each noise level according to
\eqref{eq:relative-discrepancy-score}.

For the reconstructed coefficient we use the discrete temporal norm
\begin{equation}
\|r_\delta-r\|_{L_\tau^2}
:=
\left(
\tau
\sum_{n=0}^{M-1}
\left|
r_\delta^{n+\frac12}
-
r(t_{n+\frac12})
\right|^2
\right)^{1/2}.
\label{eq:r-L2tau}
\end{equation}
The relative derivative error is defined by
\begin{equation}
E_{\mathrm{rel},z}
=
\frac{
\left(
\tau
\sum_{n=0}^{M-1}
|z^{\delta,n+\frac12}
-w'(t_{n+\frac12})|^2
\right)^{1/2}
}{
\left(
\tau
\sum_{n=0}^{M-1}
|w'(t_{n+\frac12})|^2
\right)^{1/2}
}.
\label{eq:relative-z-error}
\end{equation}

\begin{table}[ht]
\centering
\caption{Reconstruction errors for noisy measurements using adaptive
Savitzky--Golay regularised differentiation with the corrected
discrepancy criterion \eqref{eq:corrected-discrepancy-target}.}
\label{tab:noisy-reconstruction}
\begin{tabular}{ccccccc}
\hline
\(\delta\)
&
\(q\)
&
\(\eta_\delta\)
&
\(E_{\mathrm{rel},z}\)
&
\(\|r_\delta-r\|_{L^\infty}\)
&
\(\|r_\delta-r\|_{L_\tau^2}\)
&
\(\|U_\delta^M-u_h^M\|_h\)
\\
\hline
$1\%$
& 53
& $2.426\times10^{-2}$
& $1.349\times10^{-2}$
& $8.490\times10^{-3}$
& $5.068\times10^{-3}$
& $4.776\times10^{-3}$
\\
$3\%$
& 53
& $6.785\times10^{-2}$
& $3.830\times10^{-2}$
& $2.536\times10^{-2}$
& $1.419\times10^{-2}$
& $1.300\times10^{-2}$
\\
$5\%$
& 53
& $1.114\times10^{-1}$
& $6.341\times10^{-2}$
& $4.233\times10^{-2}$
& $2.348\times10^{-2}$
& $2.122\times10^{-2}$
\\
\hline
\end{tabular}
\end{table}
The results in Table~\ref{tab:noisy-reconstruction} show that the
regularised differentiation substantially reduces the instability
associated with direct differentiation of noisy measurements.
As expected, the derivative error \(\eta_\delta\) and the coefficient
reconstruction error increase as the measurement noise level
increases. Nevertheless, the reconstructed coefficient retains the
main profile of the exact function throughout the considered range of
noise levels.

The numerical state \(U_\delta^M\) is less sensitive to perturbations
than the reconstructed coefficient. This behaviour is consistent with
the structure of the inverse algorithm: the noisy data enter the
coefficient recovery through the differentiated quantity
\(z^\delta\), whereas the state is subsequently obtained by a stable
Crank--Nicolson evolution.

The results are also consistent with
Theorem~\ref{thm:inverse-convergence}, in the sense that the observed
reconstruction errors remain controlled by $\tau^2+h^2+\eta_\delta.$

This comparison should be interpreted as a numerical verification of
the conditional estimate rather than as a proof of convergence with
respect to the raw noise level \(\delta\)

The selected-window scores corresponding to the noise levels $1\%$, $3\%$, and
$5\%$ are 
$1.470\times10^{-3},
\,
1.767\times10^{-3},
\,
2.528\times10^{-3},$
respectively. These small relative scores indicate that the selected smoothing parameters closely satisfy the discrepancy criterion.

The results in Table~\ref{tab:noisy-reconstruction} show that the adaptive
regularised differentiation effectively controls the amplification of
measurement noise. For $1\%$ relative noise in $w$, the relative $L^2$ error
in the reconstructed derivative is approximately $1.35\%$, while for $3\%$
and $5\%$ measurement noise the corresponding errors are approximately
$3.83\%$ and $6.34\%$, respectively.

The resulting coefficient reconstruction errors are
\[
\|r^\delta-r\|_{L^\infty}
=
8.490\times10^{-3},
\quad
2.536\times10^{-2},
\quad
4.233\times10^{-2}
\]
for $1\%$, $3\%$, and $5\%$ noise levels, respectively. Thus, the
reconstruction error increases gradually with the level of perturbation,
while remaining moderate even for $5\%$ noise.

The numerical solution $U^\delta$ is less sensitive to the noisy
measurements. For example, for $5\%$ noise we obtain
\[
\|U_\delta^M-u_h^M\|_h
=
2.122\times10^{-2}.
\]

\begin{figure}[htbp]
\centering
\begin{subfigure}{0.48\textwidth}
\centering
\includegraphics[width=\textwidth]{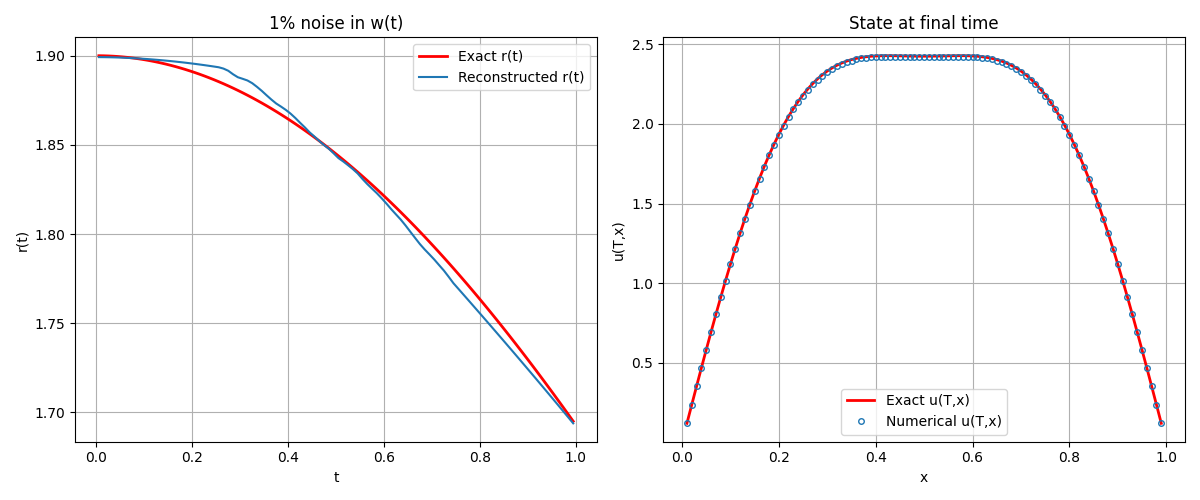}
\caption{$1\%$ noise.}
\end{subfigure}
\hfill
\begin{subfigure}{0.48\textwidth}
\centering
\includegraphics[width=\textwidth]{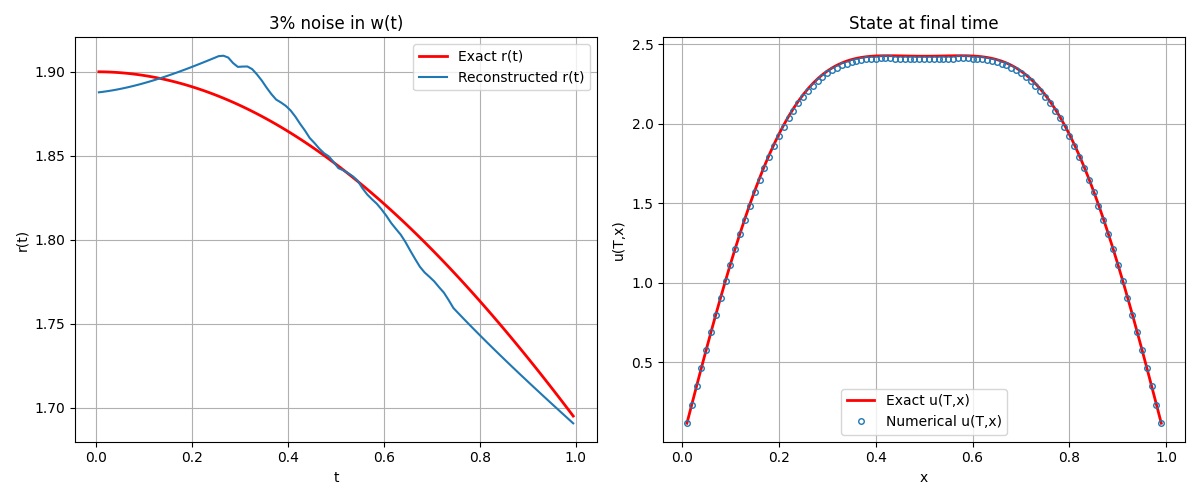}
\caption{$3\%$ noise.}
\end{subfigure}

\vspace{0.15cm}
\begin{subfigure}{0.48\textwidth}
\centering
\includegraphics[width=\textwidth]{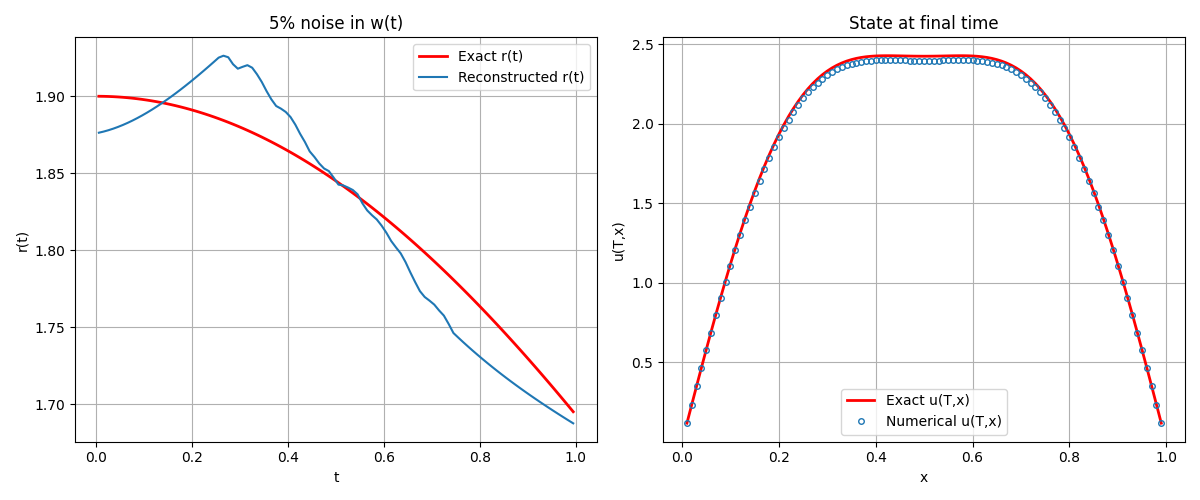}
\caption{$5\%$ noise.}
\end{subfigure}
\caption{Reconstruction of \(r(t)\) and \(u(T,x)\) from noisy integral
measurements \(w^\delta\) with \(1\%\), \(3\%\), and \(5\%\) relative
noise. The regularised derivative approximation \(z^\delta(t)\approx
w'(t)\) used in the coefficient reconstruction is obtained from
\(w^\delta\) by adaptive Savitzky--Golay differentiation.}
\label{noise1}
\end{figure}

Figure~\ref{noise1} provides a visual confirmation of the quantitative
results reported in Table~\ref{tab:noisy-reconstruction}. At the
\(1\%\) noise level, the reconstructed coefficient \(r_\delta(t)\)
closely follows the exact coefficient. As the noise level increases to
\(3\%\) and \(5\%\), the deviations become more visible, especially in
the recovered coefficient, whereas the reconstructed state
\(U_\delta^M\) remains comparatively less affected. This is consistent
with the structure of the inverse algorithm, since the perturbation in
the measured data enters the coefficient reconstruction through the
regularised derivative \(z^\delta\). The numerical behaviour is
therefore consistent with the conditional estimate
\eqref{eq:noisy-conditional-estimate} and demonstrates good empirical
robustness for the considered noise levels.

\section{Conclusion}
We studied an inverse problem for a heat equation governed by the spectral Dirichlet fractional Laplacian, in which a time-dependent source coefficient $r(t)$ is recovered from a weighted integral measurement. The spatial operator was discretised by the fractional matrix power $A_h=L_h^s$, and a Crank--Nicolson scheme was used in time. The direct scheme was shown to be unconditionally stable and second-order accurate, with error $O(\tau^2+h^2)$ under the stated spatial consistency assumptions. For the inverse problem, a closed scalar reconstruction formula for $r(t)$ was derived together with a discrete identifiability condition, and conditional convergence of the coupled state--coefficient reconstruction was established. The effect of noisy measurements was also investigated. Since direct differentiation of noisy integral data is unstable, an adaptive Savitzky--Golay regularised differentiation procedure was employed before coefficient recovery.  The noisy-data experiments demonstrate good empirical robustness and
are consistent with the conditional reconstruction estimate in terms of the regularised derivative error.

\section*{Data availability}
All data will be made available by the corresponding author upon reasonable request.

\section*{Funding}
\noindent
The research is financially supported by a grant from the Ministry of Science and Higher Education of the Republic of Kazakhstan (No. AP27508473).

\section*{Competing interests}
\noindent
The authors declare that they have no conflict of interest.

\end{document}